\documentclass[12pt]{article}
\usepackage{amsfonts,amssymb,amsmath}
\usepackage{epsfig}
\usepackage[all]{xy}
\usepackage[makeroom]{cancel}

\newcommand\RR{\mathbb R}
\newcommand\CC{\mathbb C}
\newcommand\NN{\mathbb N}
\newcommand\cT{\mathcal{T}}
\renewcommand{\Re}{\mathop{\mathrm{Re}}}
\renewcommand{\Im}{\mathop{\mathrm{Im}}}
\newcommand{\res}{\mathop{\mathrm{res}}}

\newcommand\beq{\begin{equation}}
\newcommand\eeq{\end{equation}}

\newtheorem{theorem}{Theorem}
\newtheorem{lemma}{Lemma}

\newtheorem{remark}{Remark}
\newtheorem{proposition}{Proposition}

\begin{document}

\title{Moutard transform approach to  generalized analytic functions with contour poles
\thanks{This work was started during the visits of the  first author to the Centre de Math\'ematiques 
Appliqu\'ees of \'Ecole Polytechnique in October 2015 and Institut des Hautes \'Etudes Scientifiques 
in November-December 2015. The first author was partially supported by the Russian Foundation for 
Basic Research, grant 13-01-12469 ofi-m2, by the program ``Leading scientific schools'' 
(grant NSh-4833.2014.1), by the program ``Fundamental problems of nonlinear dynamics''. }}
\author{P.G. Grinevich
\thanks{L.D. Landau Institute for Theoretical Physics,
pr. Akademika Semenova 1a, 
Chernogolovka, 142432, Russia; Lomonosov Moscow State University, Faculty of Mechanics and Mathematics, 
Russia, 119991, Moscow, GSP-1, Leninskiye Gory 1, Main Building;
Moscow Institute of Physics and Technology, 9 Institutskiy per., Dolgoprudny,
Moscow Region, 141700, Russia; e-mail: pgg@landau.ac.ru}, \and R.G. Novikov\thanks
{CNRS (UMR 7641), Centre de Math\'ematiques Appliqu\'ees, 
\'Ecole Polytechnique, 91128, Palaiseau, France; IEPT RAS, 117997, Moscow, Russia;
e-mail: novikov@cmap.polytechnique.fr}}
\date{}
\maketitle
\begin{abstract}
We continue studies of Moutard-type transforms for the generalized analytic 
functions started in \cite{GRN3}, \cite{GRN4}. In particular, we show that generalized 
analytic functions with the simplest contour poles can be Moutard transformed to the 
regular ones, at least, locally. In addition, the later Moutard-type transforms
are locally invertible. 
\end{abstract}

\section{Introduction}
We study the equations
\begin{equation}
\label{eq:ga1}
\partial_{\bar z} \psi = u \bar \psi \ \ \mbox{in} \ \ D,
\end{equation}
\begin{equation}
\label{eq:ga2}
\partial_{\bar z} \psi^+ = -\bar u \bar \psi^+ \ \ \mbox{in} \ \ D,
\end{equation}
where $D$ is an open domain in $\CC$,  $u=u(z)$ is a given function in $D$, 
$\partial_{\bar z}=\partial/\partial\bar z$.
The functions $\psi=\psi(z)$ satisfying equation (\ref{eq:ga1}) are known as generalized 
analytic functions in $D$, equation (\ref{eq:ga2}) is known as 
the conjugate equation to  (\ref{eq:ga1}); see \cite{Vek}. In the present article the 
notation $f=f(z)$ does not mean that $f$ is holomorphic. 

The classical theory of generalized analytic functions is presented in \cite{Bers}, \cite{Vek}. 
In addition, very recently in \cite{GRN3}, \cite{GRN4} it was shown that a new progress in this theory 
is possible by involving ideas of Moutard-type transforms going back to \cite{Mout}. Actually,  
ideas of Moutard-type transform were developed and successfully used in the soliton theory 
in dimension 2+1, in the spectral theory in dimension 2 and in the differential geometry; 
see \cite{DGNS}, \cite{MatvSal}, \cite{NS}, \cite{NTT}-\cite{TT} and references therein. 

We recall that in our case the Moutard-type transforms assign in quadratures to a given 
coefficient $u$ and fixed solutions $f_j$, $f^+_j$, $j=1,\ldots,N$, of equations (\ref{eq:ga1}),(\ref{eq:ga2}), 
and all formal solutions $\psi$, $\psi^+$ of (\ref{eq:ga1}), (\ref{eq:ga2}) a new coefficient $\tilde u$ and 
new related formal solutions 
$\tilde\psi$, $\tilde\psi^+$ for these generalized analytic function equations; see 
\cite{GRN3}, \cite{GRN4}. In turn, the construction of \cite{GRN3} was stimulated by recent articles by 
I.A. Taimanov \cite{Taim1}, \cite{Taim2} on the Moutard-type transforms for the Dirac operators 
in the framework of the soliton theory in dimension $2+1$.

In the classical theory of generalized analytic functions it is usually assumed that 
\begin{align}
\label{eq:cond1}
&u \in L_p(D),& \  &p>2,& \  &\mbox{if} \ \ D \ \ \mbox{is bounded},& \hspace{4cm}\\
\label{eq:cond2}
&u \in L_{p,2}(\CC),&  \ &p>2,&  \ &\mbox{if} \ \ D=\CC,& \hspace{4cm}
\end{align}
where  
\begin{align}
&L_{p,\nu}(\CC)\ \ \mbox{denotes complex-valued functions} \ \ u \ \ \mbox{such that}
\nonumber \\
&u\in  L_{p}(D_1), \ \ u_{\nu}\in  L_{p}(D_1), \ \ \mbox{where} \ \ u_{\nu}(z)=
\frac{1}{|z|^{\nu}} u\left(\frac{1}{z}\right),\\
&D_1=\{z\in\CC\, :\, |z|\le 1 \}; \nonumber
\end{align}
see \cite{Vek}.

On the other hand, one of the most important applications of the generalized analytic functions 
theory is associated with the inverse scattering in two dimensions, see \cite{BC}, 
\cite{FS}-- \cite{GRN2},  \cite{GN}, \cite{Kaz}, \cite{Rnov2}, \cite{RNov} and 
references therein.
In addition, already in \cite{GN} it was shown that, for the case of the two-dimensional 
Schr\"odinger equation, not only regular generalized analytic functions, where $u$ satisfies 
(\ref{eq:cond1}) or (\ref{eq:cond2}), but generalized analytic functions with contour poles, 
also naturally arise. However, in the latter case, assumptions  (\ref{eq:cond1}), (\ref{eq:cond2}) 
are not valid at all. It is quite likely, that the classical methods of generalized analytic 
functions dot not admit appropriate generalizations for this case. In particular,
the problem of proper solving the generalized analytic function equation 
(\ref{eq:ga1}) for $u$ with contour poles was remaining well-known, but open. 

It is in order to precise that in the framework of inverse scattering for the 
Schr\"odinger equation in two dimensions the generalized analytic function equation 
(\ref{eq:ga1}) with contour poles arises for the case when 
\begin{equation}
\label{eq:sing1}
u=\frac{u_{-1}}{\Delta}, \ \ \psi=\frac{\psi_{-1}}{\Delta}, 
\end{equation}
where $\Delta$ is real-valued and $u_{-1}$, $\psi_{-1}$, $\Delta$ are quite regular 
(e.g. real-analytic) functions on $D$, and the pole contours are the zeroes of $\Delta$.
In this setting $\psi=\psi(x,z,E)$ are the Faddeev exponentially growing solutions of 
the Schr\"odinger equation in $x\in\RR^2$ at fixed energy $E\in\RR$, $u$ is a particular 
case of the Faddeev generalized scattering data, $\Delta$ is the modified Fredholm 
determinant for the Lipman-Schwinger-Faddeev integral equation for $\psi$, and $z$ is a 
fixed-energy spectral parameter; see \cite{F1}, \cite{NKh}, \cite{GN} and references therein.

Proceeding from the aforementioned inverse scattering motivation, we consider the generalized 
analytic function equation (\ref{eq:ga1}) with contour poles for the case when $u$ is of the
form as in (\ref{eq:sing1}), and when this equation has sufficiently many (more or less as in 
the regular case) local solutions $\psi$ of the form as in (\ref{eq:sing1}). Adopting the 
terminology of \cite{GN2} we say that in this case equation (\ref{eq:ga1}) is of 
meromorphic class. Note also that equation (\ref{eq:ga1}) may be of meromorphic class only 
if principal terms of $u$ near pole contours satisfy solvability conditions; for simple 
contour poles these conditions were found in \cite{GN}.

Actually, our recent works \cite{GRN3}, \cite{GRN4} were motivated considerably by the aforementioned open 
problem of proper solving the generalized analytic function equation (\ref{eq:ga1}) with contour
poles for the meromorphic case. In particular, in \cite{GRN3} we give examples of coefficients $u$
when equation (\ref{eq:ga1}) is of meromorphic class and can be efficiently studied 
using Moutard-type transforms. In addition, we were stimulated by results of  \cite{APP}, \cite{Crum},
\cite{GN1}, \cite{NT}, \cite{TT}, \cite{Taim1} on efficient 
applications of Darboux-Crum and Moutard-type transforms to studies of some important linear 
ODE's and PDS's with singular coefficients.

The results of the present work can be summarized as follows: 

\begin{itemize}
\item We give composition and inversion formulas for the simple Moutard-type transforms for the 
equations (\ref{eq:ga1}), (\ref{eq:ga2}); see Theorems~\ref{th:1},~\ref{th:2} in Subsection~\ref{sec:2.2};
\item We show that any equation (\ref{eq:ga1}) of meromorphic class with a simple contour pole 
can be transformed to a regular one in a neighborhood of the pole contour via an appropriate
simple Moutard-type transform; see Theorem~\ref{th:3} in Subsection~\ref{sec:3.4}.
\end{itemize}

\section{Moutard-type transforms for generalized analytic functions}
\label{sec:2}
\subsection{Basic construction}
Let 
\begin{equation}
\label{eq:sol1}
f_j=f_j(z)\ \ \mbox{and} \ \ f^+_j=f^+_j(z), \ \ j=1,\ldots,N,
\end{equation}
denote a set of fixed solutions of equations (\ref{eq:ga1}) and (\ref{eq:ga2}), respectively. 

Let $\psi$, $\psi^+$ be arbitrary solutions of  (\ref{eq:ga1}), (\ref{eq:ga2}). We define  
an imaginary-valued potential $\omega_{\psi,\psi^+}=\omega_{\psi,\psi^+}(z)$ such that
\begin{equation}
\label{eq:kern1}
\partial_{z} \omega_{\psi,\psi+} =\psi\psi^+, \ \ 
\partial_{\bar z} \omega_{\psi,\psi^+} =-\overline{\psi\psi^+} \ \ \mbox{in} \ \ 
D.
\end{equation}
We recall that this definition is self-consistent, at least, under the assumption that 
$D$ is simply-connected. The integration constant is imaginary-valued and may depend on concrete situation. 

As in \cite{GRN3}, \cite{GRN4} we consider imaginary-valued potentials $\omega_{f_j,f^+_k}$, 
$j,k=1,\ldots,N$, and we set
\begin{equation}
\label{eq:omega1}
\Omega=\begin{bmatrix} 
\omega_{f_1,f^+_1} & \omega_{f_2,f^+_1} & \ldots & \omega_{f_N,f^+_1}  \\
\omega_{f_1,f^+_2} & \omega_{f_2,f^+_2} & \ldots & \omega_{f_N,f^+_2}  \\
\vdots & \vdots & \ddots & \vdots \\
\omega_{f_1,f^+_N} & \omega_{f_2,f^+_N} & \ldots & \omega_{f_N,f^+_N}  \\
\end{bmatrix}.
\end{equation}
Following \cite{GRN3}, for equations (\ref{eq:ga1}), (\ref{eq:ga2}) we consider the Moutard-type transform 
\begin{align}
&\xymatrix{{\mathstrut u}\ar[r]^-{\mathcal{M}} & {\mathstrut\tilde{u}}}, \label{eq:moutard1} \\ 
&\xymatrix{{\mathstrut \{\psi,\psi^+\}}\ar[r]^-{\mathcal{M}} & {\mathstrut\{\tilde{\psi},\tilde{\psi^+}\} }},  
\label{eq:moutard2}\\
&\mathcal{M}=\mathcal{M}_{u,f,f^+}= \mathcal{M}_{u,f_1,\ldots,f_N,f_1^+,\ldots,f_N^+}  \label{eq:moutard2.5}\
\end{align}
defined as follows:
\begin{equation}
\label{eq:moutard3}
\tilde u = u + \begin{bmatrix}f_1 & \ldots & f_N \end{bmatrix} 
 \Omega^{-1} \begin{bmatrix} \overline{f_1^+} \\ \vdots \\ 
\overline{f_N^+}  \end{bmatrix},
\end{equation}
\begin{equation}
\label{eq:moutard5}
\tilde\psi =
\psi -
\begin{bmatrix}f_1 & \ldots & f_N \end{bmatrix} 
\Omega^{-1}
\begin{bmatrix}\omega_{\psi,f^+_1} \\ \vdots \\ \omega_{\psi,f^+_N} 
\end{bmatrix},
\end{equation}
\begin{equation}
\label{eq:moutard6}
\tilde\psi^+= \psi^+ -
\begin{bmatrix}f_1^+ & \ldots & f_N^+ \end{bmatrix} 
(\Omega^{-1})^t 
\begin{bmatrix}\omega_{f_1,\psi^+} \\ \vdots \\ \omega_{f_N,\psi^+} 
\end{bmatrix},
\end{equation}
where $\psi$,  $\psi^+$ are formal solutions to equations 
(\ref{eq:ga1}), (\ref{eq:ga2}), respectively, $\omega_{\psi,f^+_j}$ and $\omega_{f_j,\psi^+}$ are defined as 
in (\ref{eq:kern1}); $t$ in (\ref{eq:moutard6}) stands for the matrix transposition.

Due to results of \cite{GRN3}, the transformed functions ${\tilde\psi}$,  ${\tilde\psi}^+$ 
solve the transformed generalized-analytic function equations:
\begin{align}
\label{eq:ga3}
&\partial_{\bar z} \tilde\psi= \tilde u\, \overline{\psi}& \ \ &\mbox{in} \ \ D,\\
\label{eq:ga4}
&\partial_{\bar z} \tilde\psi^+= -\overline{\tilde u}\, 
\overline{\tilde\psi^+}& \ \ 
&\mbox{in} \ \ D.
\end{align}
\subsection{Composition and inversion of simple Moutard  transforms}
\label{sec:2.2}

We say that the Moutard-type transforms (\ref{eq:moutard1})- (\ref{eq:moutard6}) are simple 
if $N=1$. In this case:
\begin{equation}
\label{eq:s2:moutard1}
\tilde\psi=
\psi-
f_1 \frac{\omega_{\psi,f_1^+}}{\omega_{f_1,f_1^+}},
\end{equation}
\begin{equation}
\label{eq:s2:moutard2}
\tilde{\psi}^+= \psi^+ - f_1^+  \frac{\omega_{f_1,\psi^+}}{\omega_{f_1,f_1^+}},
\end{equation}
\begin{equation}
\label{eq:s2:moutard3}
\tilde u = u + \frac{f_1\overline{f_1^+}}{\omega_{f_1,f_1^+}}.
\end{equation}
\begin{proposition}
\label{pr:1}
For a simple Moutard transform (\ref{eq:s2:moutard1})- (\ref{eq:s2:moutard3}) the following formula holds:
\beq
\label{eq:trans:pot}
\omega_{\tilde\psi,\tilde\psi^+}=\frac{\omega_{\psi,\psi^+}\omega_{f_1,f_1^+}-\omega_{\psi,f_1^+}
\omega_{f_1,\psi^+}}{\omega_{f_1,f_1^+}}+c_{\tilde\psi,\tilde\psi^+},
\eeq 
where  $\omega_{\tilde\psi,\tilde\psi^+}$ is defined according to (\ref{eq:kern1}), 
and $c_{\tilde\psi,\tilde\psi^+}$ is an imaginary-valued integration constant.
\end{proposition}

Proposition~\ref{pr:1} was announced in \cite{GRN4} and is proved in Section~\ref{sec:4} of the present work. 

Let $f_1$, $f_2$ and $f_1^+$, $f_2^+$ be some fixed solutions of equations 
(\ref{eq:ga1}) and (\ref{eq:ga2}), respectively, with given $u$. Let 
\begin{align}
&\xymatrix{{\mathstrut u}\ar[r]^-{\mathcal{M}} & {\mathstrut\tilde{\tilde u}}} ,  
\label{eq:pr2:moutard1} \\ 
&\xymatrix{ \left\{\vphantom{{\tilde{\tilde\psi}} } \psi,\psi^+\right\} 
\ar[r]^-{\mathcal{M}} & \left\{ {\tilde{\tilde\psi}}, 
{\tilde{\tilde\psi}}^+ \right\} }\label{eq:pr2:moutard2},
\end{align}

\begin{equation}
\label{eq:comp1}
\mathcal{M} = \mathcal{M}_2 \circ \mathcal{M}_1,
\end{equation}
where
\begin{enumerate}
\item
$\mathcal{M}_1$ is the simple Moutard transform for equations
(\ref{eq:ga1}), (\ref{eq:ga2}) with coefficient $u$,  generated by 
$f_1$, $f_1^+$ and given by formulas (\ref{eq:s2:moutard1})-(\ref{eq:s2:moutard3});
\item
$\mathcal{M}_2$ is the simple Moutard transform for equations
(\ref{eq:ga3}), (\ref{eq:ga4}) with coefficient $\tilde u =\mathcal{M}_1 u$, given by:
\begin{equation}
\label{eq:pr2:moutard3}
\tilde{\tilde\psi}=
\tilde{\psi}-
\tilde{f}_2\ \frac{\omega_{\tilde\psi,\tilde{f}_2^+}}{\omega_{\tilde{f}_2,\tilde{f}_2^+}},
\end{equation}
\begin{equation}
\label{eq:pr2:moutard4}
\tilde{\tilde\psi}^+= \tilde\psi^+ - \tilde{f}_2^+\  
\frac{\omega_{\tilde{f}_2,\tilde\psi^+}}{\omega_{\tilde{f_2},\tilde{f}_2^+}},
\end{equation}
\begin{equation}
\label{eq:pr2:moutard5}
\tilde{\tilde u} = \tilde u + \frac{\tilde f_2\overline{\tilde f_2^+}}
{\omega_{\tilde{f}_2,\tilde{f}_2^+}},
\end{equation}
where we assume that: 
\begin{itemize}
\item ${\tilde\psi}=\mathcal{M}_1\psi$, ${\tilde\psi}^+=\mathcal{M}_1\psi^+$, 
${\tilde f}_2=\mathcal{M}_1f_2$, ${\tilde f}_2^+=\mathcal{M}_1f_2^+$;
\item $\omega_{\tilde\psi,\tilde{f}_2^+}$, $\omega_{\tilde{f}_2,\tilde\psi^+}$ and  $\omega_{\tilde{f}_2,\tilde{f}_2^+}$       
are given by (\ref{eq:trans:pot}) with $c_{\tilde\psi,\tilde{f}_2^+}$, $c_{\tilde{f}_2,\tilde\psi^+}$ and
$c_{\tilde{f}_2,\tilde{f}_2^+}$ equal to zero. 
\end{itemize}
\item We assume that $\omega_{f_1,f_1^+}\ne 0$,  $\omega_{\tilde{f}_2,\tilde{f}_2^+}\ne 0$. 
\end{enumerate}
\begin{theorem}
\label{th:1}
Let $f_1$, $f_2$ and $f_1^+$, $f_2^+$ be some 
fixed solutions of equations (\ref{eq:ga1}) and (\ref{eq:ga2}), respectively, with given $u$,
and let $\mathcal{M}$ be defined as in (\ref{eq:pr2:moutard1})-(\ref{eq:comp1}).
Then  $\mathcal{M}$ coincides with the  Moutard transform given by formulas 
(\ref{eq:moutard1})-(\ref{eq:moutard6}) for $N=2$ and generated by the initial functions 
$f_1$, $f_2$, $f_1^+$, $f_2^+$.
\end{theorem}
Schematically, the result of Theorem~\ref{th:1} can be also presented as follows:
\begin{equation}
\label{eq:pr2:1}
\xymatrix{ \psi
\ar[rrr]^-{\{f_1,f_1^+ \}}_{\mathcal{M}_1} 
\ar@/_2pc/[rrrrrr]^-{\vphantom{\big(}
\{f_1,f_2, f_1^+, f_2^+ \}}
  & & &
\tilde\psi
\ar[rrr]^-{\{{\tilde{f}}_2,{\tilde{f}}_2^+\}}_{\mathcal{M}_2}
& & &  
\tilde{\tilde\psi}
}
\end{equation}

Theorem~\ref{th:1} is proved in Section~\ref{sec:4}. 

Next, in order to inverse simple Moutard transforms, we consider:
\begin{equation}
\label{eq:s3:moutard4}
f_2\equiv 0,\ \  f_2^+\equiv 0, \ \omega_{f_2,f_1^+} = i, \ \ \omega_{f_1,f_2^+} = i;
\end{equation}
\begin{equation}
\label{eq:s3:moutard4.5}
\tilde f_2=\mathcal{M}_1  f_2 ={\hat{f}}= -i\frac{f_1 }{\omega_{f_1,f_1^+}}, 
\ \ 
\tilde f_2^+=\mathcal{M}_1 f_2^+={\hat{f}}^+= -i \frac{f_1^+}{\omega_{f_1,f_1^+}}.
\end{equation}
In particular, ${\hat{f}}$ and ${\hat{f}}^+$ defined in (\ref{eq:s3:moutard4.5}) satisfy equations 
(\ref{eq:ga3}) and (\ref{eq:ga4}), respectively,  with the coefficient $\tilde u$  given by (\ref{eq:s2:moutard3}).

\begin{theorem}
\label{th:2}
Let $\mathcal{M}_1$ be a simple Moutard transform defined as in 
(\ref{eq:s2:moutard1})-(\ref{eq:s2:moutard3}), where $\omega_{f_1,f_1^+}\ne 0$. Let $\mathcal{M}_2$ be the 
simple Moutard transform for equations (\ref{eq:ga3}), (\ref{eq:ga4}) with coefficient 
$\tilde u=\mathcal{M}_1 u$, given by (\ref{eq:pr2:moutard3})-(\ref{eq:pr2:moutard5}), where
${\tilde\psi}=\mathcal{M}_1\psi$, ${\tilde\psi}^+=\mathcal{M}_1\psi^+$, 
${\tilde{f}}_2=\hat{f}$, ${\tilde{f}}_2^+=\hat{f}^+$, where  ${\hat{f}}$ and  ${\hat{f}}^+$ 
are defined in (\ref{eq:s3:moutard4.5}).
Then:
\begin{enumerate}
\item The potentials $\omega_{\tilde\psi,\tilde{f}_2^+}$,  $\omega_{\tilde{f}_2,\tilde\psi^+}$, 
$\omega_{\tilde{f}_2,\tilde{f}_2^+}$ can be chosen as follows:
\begin{equation}
\label{eq:s3:moutard5}
\omega_{\tilde\psi,\tilde{f}_2^+}=-i\frac{\omega_{\psi,f_1^+}}{\omega_{f_1,f_1^+}}, \ \ 
\omega_{\tilde{f}_2,\tilde\psi^+}=-i\frac{\omega_{f_1,\psi^+}}{\omega_{f_1,f_1^+}}, \ \ 
\omega_{\tilde{f}_2,\tilde{f}_2^+}=\frac{1}{\omega_{f_1,f_1^+}};
\end{equation}
\item Under assumptions (\ref{eq:s3:moutard5}), the composition $\mathcal{M}=\mathcal{M}_2  \circ \mathcal{M}_1$ is the 
identical transformation.
\end{enumerate}
\end{theorem}

Thereom~\ref{th:2} is proved in Section~\ref{sec:4}.

\section{Removing the simplest contour pole singularity}
\subsection{Real analytic pole contour}
We consider a real analytic curve $\Gamma\subset D$:
\begin{equation}
\label{eq:curve1}
\Gamma=\{z(\tau):\tau\in]it_1,it_2[ \},  \ \ t_1,t_2\in\RR,
\end{equation}
where:
\begin{itemize}
\item $z(\tau)\in D$ for $\tau\in]it_1,it_2[$,
\item $z$ is (complex-valued) real-analytic on $]it_1,it_2[$,
\item $z(\tau'_1)\ne z(\tau'_2)$ for $\tau'_1\ne\tau'_2$,
\item $d z(it)/dt\ne 0$ for for $t\in]t_1,t_2[$.
\end{itemize}
Here $D$ is the domain in (\ref{eq:ga1}), (\ref{eq:ga2}). 

Let 
\begin{equation}
\label{eq:curve2}
\cT=\cT_{a,b,\varepsilon}=\{\tau\in\CC: \ a<\Im\tau<b, \ \left|\Re\tau\right|<\varepsilon \},
\end{equation}
where $a,b,\varepsilon\in\RR$, $\varepsilon>0$. 

As a corollary of our assumptions, $z$ in (\ref{eq:curve1}) could be continued to a holomorphic bijection 
$Z$: 
\begin{align}
\label{eq:curve3}
& Z:\ \cT_{a,b,\varepsilon} \rightarrow D_{a,b,\varepsilon}, \ \ \tau\rightarrow z(\tau), \\
& Z^{-1}:\ D_{a,b,\varepsilon} \rightarrow \cT_{a,b,\varepsilon}, \ \ z\rightarrow \tau(z),\nonumber 
\end{align}
for some $a,b,\varepsilon$ such that $t_1<a<b<t_2$, \ $\varepsilon>0$, where $D_{a,b,\varepsilon}\subset D$.

Actually, we consider $\Gamma$ of (\ref{eq:curve1}) as a pole contour for equations (\ref{eq:ga1}), (\ref{eq:ga2}).

\subsection{Holomorphic change of variables}

In the domain $ D_{a,b,\varepsilon}$ of (\ref{eq:curve3}) we rewrite equations  (\ref{eq:ga1}), (\ref{eq:ga2}) in 
variable $\tau\in\cT_{a,b,\varepsilon}$. In this connection, following \cite{GRN4} we consider
\begin{align}
\label{eq:conf2}
&u_*(\tau)=u(z(\tau))\,\left|\frac{\partial z(\tau)}{\partial \tau}\right|, \\ 
&\psi_*(\tau)=\psi(z(\tau))\,\sqrt{\frac{\partial z(\tau)}{\partial \tau}},\ \ \label{eq:conf3}
\psi^+_*(\tau)=\psi^+(z(\tau))\,\sqrt{\frac{\partial z(\tau)}{\partial \tau}}, 
\end{align}
where $u(z)$, $\psi(z)$,  $\psi^+(z)$ are the same that in equations (\ref{eq:ga1}), 
(\ref{eq:ga2}), $z(\tau)$ is the same that in (\ref{eq:curve3}), $\tau\in\cT_{a,b,\varepsilon}$. Then 
(see \cite{GRN4}):
\begin{align}
&\partial_{\bar\tau} \psi_* = u_* \bar \psi_* \ \ \mbox{in} \ \  \cT_{a,b,\varepsilon}, \label{eq:conf4}\\
&\partial_{\bar\tau} \psi^+_* = -\bar u_* \bar \psi^+_* \ \ \mbox{in} \ \ \cT_{a,b,\varepsilon}. \label{eq:conf5}
\end{align}

In addition, we have (see \cite{GRN4}):
\begin{equation}
\label{eq:conf2.5}
\mathcal{M}_{u_*,f_*,f^+_*} \circ Z^{-1} = Z^{-1} \circ \mathcal{M}_{u,f,f^+},
\end{equation}
where: 
\begin{itemize} 
\item $Z^{-1}$ is considered as a map of the conjugate pairs of equations  (\ref{eq:ga1}), (\ref{eq:ga2}) in 
$D_{a,b,\varepsilon}$ into the conjugate pairs of equations  (\ref{eq:conf4}), (\ref{eq:conf5})  in 
$\cT_{a,b,\varepsilon}$ and is defined according to  (\ref{eq:conf2}), (\ref{eq:conf3});
\item $\mathcal{M}_{u,f,f^+}$ for (\ref{eq:ga1}), (\ref{eq:ga2}) in $D_{a,b,\varepsilon}$ and 
$\mathcal{M}_{u_*,f_*,f^+_*}$ for (\ref{eq:conf4}), (\ref{eq:conf5})  in $\cT_{a,b,\varepsilon}$
are defined according to formulas (\ref{eq:moutard1})- (\ref{eq:moutard6}), where $u_*$,
$f_*=\{f_{1,*},\ldots,f_{N,*} \}$, $f^+_*=\{f_{1,*}^+,\ldots,f_{N,*}^+ \}$ are defined according to 
(\ref{eq:conf2}), (\ref{eq:conf3}), and 
\begin{equation}
\label{eq:conf2.6}
\omega_{\psi_*,\psi^+_*}(\tau)= \omega_{\psi,\psi^+}(z(\tau)),
\end{equation}
for all involved potentials.
\end{itemize}
In the framework of the Moutard transform approach, using the commutativity relation 
(\ref{eq:conf2.5}) we reduce local studies of equations (\ref{eq:ga1}), (\ref{eq:ga2}) with 
contour pole at $\Gamma$ in (\ref{eq:curve1}) to the case of contour pole at the straight line 
\begin{equation}
\label{eq:conf2.7}
\Gamma_*=\{\tau\in \cT_{a,b,\varepsilon}: \ \Re\tau=0 \},
\end{equation}
where $\cT_{a,b,\varepsilon}$ is defined as in (\ref{eq:curve2}), (\ref{eq:curve3}).
\begin{remark} We recall that, in view of formulas (\ref{eq:conf2}),  
(\ref{eq:conf3}), the generalized analytic functions  $\psi$,  $\psi^+$ 
of  (\ref{eq:ga1}), (\ref{eq:ga2})  can be treated as spinors, i.e. 
differential forms of type $\left(\frac{1}{2},0\right)$, and $u$ can be 
treated as differential form of type  $\left(\frac{1}{2},\frac{1}{2}\right)$;  
see \cite{GRN4}. These forms can be written as:
\begin{equation}
u=u(z)\sqrt{dzd\bar z}, \ \ \psi=\psi(z)\sqrt{dz}, \ \  \psi^+=\psi^+(z)\sqrt{dz}.
\end{equation}
\end{remark}

\subsection{Constraints on the meromorphic class 
coefficients at singularity}
\label{sec3.3}

We consider equations (\ref{eq:conf4}), (\ref{eq:conf5}) in 
$\cT=\cT_{a,b,\varepsilon}$ defined by (\ref{eq:curve2}) for the case of simplest 
pole at $\Gamma_*$ defined by (\ref{eq:conf2.7}). We write $\tau=x+iy$, 
$\bar\tau=x-iy$.

We assume that 
\begin{equation}
\label{eq:sing2}
u_*(\tau)=e^{2i\phi(y)}\sum\limits_{j=-n}^{+\infty}
r_j(y)\, x^j\ \ \mbox{in}\ \  \cT, 
\end{equation}
where $\phi$, $r_j$ are quite regular functions on the interval 
$]a,b[\approx\Gamma_*$, and $\phi$, $r_{-n}$ are real-valued, $n\in\NN$. For 
this case we consider solutions $\psi_*$, $\psi_*^+$ of (\ref{eq:conf4}), 
(\ref{eq:conf5}) near $\Gamma_*$ of the following form:
\begin{equation}
\label{eq:sing3}
\psi_*=\sum\limits_{j=-n'}^{+\infty} \alpha_j(y)\, x^j, \ \ 
\psi_*^+=\sum\limits_{j=-n''}^{+\infty} \alpha_j^+(y)\, x^j, \ \ n',n''\in\NN,
\end{equation}
where $\alpha_j$, $\alpha_j^+$ are quite regular on $]a,b[$ and $\alpha_{-n'}$, $\alpha_{-n''}^+$ 
are not-zero almost everywhere at $]a,b[$.

Note that in this section we consider $u_*$, $\psi_*$, $\psi_*^+$ in formulas 
(\ref{eq:sing2}), (\ref{eq:sing3}) as formal power series in variable $x$.  

\begin{lemma}
\label{lem:1}
Assume that equation (\ref{eq:conf4}) with coefficient $u_*$ as in 
(\ref{eq:sing2}) has, at least, one solution $\psi_*$ of the form as in 
(\ref{eq:sing3}). Then  $n=1$ and $|r_{-1}(y)|\equiv n'/2$ in (\ref{eq:sing2}).
\end{lemma}

Lemma~\ref{lem:1} follows from formal substitutions of (\ref{eq:sing2}), 
(\ref{eq:sing3}) into  (\ref{eq:conf4}),  (\ref{eq:conf5}).  

In the present article we restrict ourselves to the simplest (but generic) 
case when $n'=n''=1$. In this case, without loss of generality, 
we can assume that $r_{-1}=-1/2$, adding $\pi/2$ to the phase $\phi$, if 
necessary.
%\begin{equation}
%\label{eq:sing4}
%u_*(\tau)=e^{2i\phi(y)}\left(-\frac{1}{2x}+r_0(y)+ r_1(y) x+  r_2(y) x^2+
%\ldots  \right)  
%\end{equation}

Then, as a corollary of results of \cite{GN}, equation 
(\ref{eq:conf4}) is of meromorphic class, at least formally, if and only if:
\begin{equation}
\label{eq:sing5}
\Re r_0(y) \equiv 0, \ \ y\in]a,b[,
\end{equation}
\begin{equation}
\label{eq:sing6}
\Im r_1(y)=\frac{1}{2} \frac{d^2 \phi(y)}{dy^2}, \ \ y\in]a,b[.
\end{equation}
For completeness of exposition, this result is proved in Section~\ref{sec:5}.

Here, belonging of equation (\ref{eq:conf4}) to meromorphic class means that  
equation (\ref{eq:conf4}) has local solutions $\psi_*$ near $\Gamma_*\subset
\cT$ of the form (\ref{eq:sing3}) with $n'=1$ parametrised by two real-values 
functions (one complex-valued function) on $\Gamma_*$, i.e., roughly speaking,
equation (\ref{eq:conf4}) has as many local solutions $\psi_*$ near $\Gamma_*$
as in the regular case.

Actually, equation (\ref{eq:conf4}) and formulas (\ref{eq:sing2}), 
(\ref{eq:sing3}) for $u_*$, $\psi_*$ with $n=1$, $n'=1$, $r_{-1}=-1/2$ imply 
that 
\begin{equation}
\label{eq:sing7}
\Im e^{-i\phi(y)}\alpha_{-1}(y) = 0, \ \ y\in]a,b[.
\end{equation}
In addition, under conditions (\ref{eq:sing5}), (\ref{eq:sing6}), the solutions 
$\psi_*$ of  (\ref{eq:conf4}), (\ref{eq:sing3}) are parametrised by 
$\beta_{-1}(y)$ and $\Im \beta_1(y)$ on $]a,b[$, where 
$\alpha_j(y)=e^{i\phi(y)}\beta_j(y)$; see Section~\ref{sec:5}.

Finally, one can see that the meromorphic class conditions (\ref{eq:sing5}), 
(\ref{eq:sing6}) for equation (\ref{eq:conf4}) imply the related meromorphic 
class conditions for the conjugate equation (\ref{eq:conf5}).

\subsection{Moutard transform to the regular case}
\label{sec:3.4}

We consider equations  (\ref{eq:conf4}),  (\ref{eq:conf5}) for the case when 
\begin{equation}
\label{eq:sing7.5}
u_*=e^{2i\phi(y)}\left(-\frac{1}{2x}+r_0(y)+r_1(y) x+ O\left(x^2\right) \right) \ \ \mbox{in} \ \ \cT\cup\partial\cT,
\end{equation}
where
\begin{align}
\nonumber
&\phi \ \ \mbox{is real-valued}, \ \ \phi\in C^2\left([a,b]\right),\\
\label{eq:sing8}
& u_* +\frac{e^{2i\phi(y)}}{2x}\in    C^1\left(\cT\cup\partial\cT\right),\\
& r_0, r_1 \ \ \mbox{satisfy  (\ref{eq:sing5}), (\ref{eq:sing6})},\nonumber
\end{align}
$\cT$ is defined by (\ref{eq:curve2}).

We assume that equations  (\ref{eq:conf4}),  (\ref{eq:conf5}) have some 
solutions $f_*$, $f_*^+$ such that 
\begin{align}
\nonumber
&f_*=e^{i\phi(y)}\left(\frac{\beta_{-1}(y)}{x}+\beta_0(y)+O(x) \right),\\
&\beta_{-1}\ \ \mbox{is real-valued}, \ \ 
\beta_{-1}>0,  \ \ 
\beta_{-1}\in C^1\left([a,b]\right),\label{eq:sing9}\\
& f_* -\frac{e^{i\phi(y)}\beta_{-1}(y)}{x}\in  C^1\left(\cT\cup\partial\cT\right),\nonumber
\end{align}
\begin{align}
\nonumber
&f_*^+=e^{-i(\phi(y)+\pi/2)}\left(\frac{\beta^+_{-1}(y)}{x}+\beta^+_{0}(y)+O(x) 
\right),\\
&\beta^+_{-1}\ \ \mbox{is real-valued}, \ \ 
\beta^+_{-1}>0,  \ \ 
\beta^+_{-1}\in C^1\left([a,b]\right),\label{eq:sing10}\\
& f^+_* -\frac{e^{-i(\phi(y)+\pi/2)}\beta^+_{-1}(y)}{x}\in  C^1\left(\cT\cup\partial\cT\right),\nonumber
\end{align}
where $\partial_x O(x)=O(1)$,  $\partial_y O(x)=O(x)$. 

Note that from point of view of formal considerations of 
Subsection~\ref{sec3.3} such solutions $f_*$, $f_*^+$ always exist. 

\begin{theorem}
\label{th:3}
Let $u_*$ satisfy (\ref{eq:sing7.5}), (\ref{eq:sing8}) and equations 
(\ref{eq:conf4}),  (\ref{eq:conf5}) have some solutions $f_*$, $f_*^+$ 
satisfying (\ref{eq:sing9}), (\ref{eq:sing10}). Let $\omega_{f_*,f_*^+}$ 
be some potential defined according to (\ref{eq:kern1}). Then 
\begin{equation}
\label{eq:th3:1}
\tilde u_* = \mathcal{M}_{u_*,f_*,f_*^+} u_*=O(1) \ \ \mbox{in} \ \ \cT_{a,b,\delta}\cup\partial\cT_{a,b,\delta}
\end{equation}
for some $\delta\in]0,\varepsilon[$, where  $\mathcal{M}_{u_*,f_*,f_*^+}$ is defined according to formulas 
(\ref{eq:moutard1})-(\ref{eq:moutard6}) for $N=1$.

Theorem~\ref{th:3} is proved in Section~\ref{sec:6}.
\end{theorem}

The point is that that the results of Theorems~\ref{th:2},~\ref{th:3} and the commutativity formula 
(\ref{eq:conf2.5}) reduce local studies of equations 
(\ref{eq:ga1}), (\ref{eq:ga2}) near the simplest contour pole singularity to the regular case. 

\section{Proofs of Proposition~\ref{pr:1} and Thereoms~\ref{th:1}-\ref{th:2}}
\label{sec:4}

\subsection{Proof of  Proposition~\ref{pr:1}}
Let $\omega_{\tilde\psi,\tilde\psi^+}$ be given by (\ref{eq:trans:pot}). Then
it is sufficient to show that
\begin{equation}
\label{eq:ppr1:1}
\partial_{z} \omega_{\tilde\psi,\tilde\psi^+} =
\tilde\psi\tilde\psi^+, \ \ 
\partial_{\bar z} \omega_{\tilde\psi,\tilde\psi^+} =
-\overline{\tilde\psi}\overline{\tilde\psi^+}.
\end{equation}
Using the definitions of $\omega_{\psi,\psi^+}$,  $\omega_{\psi,f_1^+}$,  
$\omega_{f_1,\psi^+}$,  $\omega_{f_1,f_1^+}$, and ${\tilde\psi}$,  ${\tilde\psi}^+$  we have:
$$
\partial_{z} \omega_{\tilde\psi,\tilde\psi^+} =
\partial_z\left(\frac{\omega_{\psi,\psi^+}\omega_{f_1,f_1^+}-\omega_{\psi,f_1^+}\omega_{f_1,\psi^+}}{\omega_{f_1,f_1^+}}+c_{\tilde\psi,\tilde\psi^+}\right)=
$$
$$
=\partial_z\omega_{\psi,\psi^+}-\partial_z\omega_{\psi,f_1^+} \frac{\omega_{f_1,\psi^+}}{\omega_{f_1,f_1^+}}-\partial_z\omega_{f_1,\psi^+} \frac{\omega_{\psi,f_1^+}}{\omega_{f_1,f_1^+}}+
\partial_{z}\omega_{f_1,f_1^+} \frac{\omega_{f_1,\psi^+}}{\omega_{f_1,f_1^+}}\frac{\omega_{\psi,f_1^+}}{\omega_{f_1,f_1^+}}=
$$
$$
=\psi\psi^+ - \psi f_1^+ \frac{\omega_{f_1,\psi^+}}{\omega_{f_1,f_1^+}}-
f_1\psi^+  \frac{\omega_{\psi,f_1^+}}{\omega_{f_1,f_1^+}}+
\psi\psi^+  \frac{\omega_{f_1,\psi^+}}{\omega_{f_1,f_1^+}}\frac{\omega_{\psi,f_1^+}}{\omega_{f_1,f_1^+}}=
$$
$$
=\left(\psi - f_1 \frac{\omega_{\psi,f_1^+}}{\omega_{f_1,f_1^+}}  \right) \left(\psi^+ -f_1^+ 
\frac{\omega_{f_1,\psi^+}}{\omega_{f_1,f_1^+}} \right)=
\tilde\psi \tilde\psi^+;
$$
$$
\partial_{\bar z} \omega_{\tilde\psi,\tilde\psi^+} =
\partial_{\bar z}\left(\frac{\omega_{\psi,\psi^+}\omega_{f_1,f_1^+}-
\omega_{\psi,f_1^+}\omega_{f_1,\psi^+}}{\omega_{f_1,f_1^+}}+c_{\tilde\psi,\tilde\psi^+}\right)=
$$
$$
=\partial_{\bar z}\omega_{\psi,\psi^+}-\partial_{\bar z}\omega_{\psi,f_1^+} 
\frac{\omega_{f_1,\psi^+}}{\omega_{f_1,f_1^+}}-\partial_{\bar z}\omega_{f_1,\psi^+} \frac{\omega_{\psi,f_1^+}}{\omega_{f_1,f_1^+}}+
\partial_{z}\omega_{f_1,f_1^+} \frac{\omega_{f_1,\psi^+}}{\omega_{f_1,f_1^+}}\frac{\omega_{\psi,f_1^+}}{\omega_{f_1,f_1^+}}=
$$
$$
=-\overline{\psi}\overline{\psi^+} + \overline{\psi}\overline{f_1^+} 
\frac{\omega_{f_1,\psi^+}}{\omega_{f_1,f_1^+}}+ \overline{f_1}\overline{\psi^+}  
\frac{\omega_{\psi,f_1^+}}{\omega_{f_1,f_1^+}}- \overline{\psi}\overline{\psi^+}  
\frac{\omega_{f_1,\psi^+}}{\omega_{f_1,f_1^+}}\frac{\omega_{\psi,f_1^+}}{\omega_{f_1,f_1^+}}=
$$
$$
=-\left(\overline{\psi} - \overline{f_1} \frac{\omega_{\psi,f_1^+}}{\omega_{f_1,f_1^+}}  \right) \left(\overline{\psi^+} -\overline{f_1^+} \frac{\omega_{f_1,\psi^+}}{\omega_{f_1,f_1^+}} \right)=
-\overline{\tilde\psi} \overline{\tilde\psi^+}.
$$
Thus, the proof of Proposition~\ref{pr:1} is completed.

\subsection{Proof of  Theorem~\ref{th:1}}
Due to formula (\ref{eq:trans:pot}) for $\omega_{\tilde f_2,\tilde f_2^+}$ with 
$c_{\tilde f_2,\tilde f_2^+}=0$ and due to the assumptions that $\omega_{f_1,f_1^+}\ne0$,
$\omega_{\tilde f_2,\tilde f_2^+}\ne0$, we have: 
\begin{equation}
\label{eq:ppr2:1}
\det\Omega=\omega_{f_1,f_1^+}\omega_{\tilde f_2,\tilde f_2^+}\ne0,
\end{equation}
where $\Omega$ is defined according to (\ref{eq:omega1}) for $N=2$.

Due to the definition of $\mathcal{M}$ according to (\ref{eq:pr2:moutard1})-(\ref{eq:comp1}),
using (\ref{eq:trans:pot}) for 
$\omega_{\tilde\psi,\tilde\psi^+}=\omega_{\tilde f_2,\tilde f_2^+}$, $\omega_{\tilde\psi,\tilde f_2^+}$ 
with $c_{\tilde\psi,\tilde\psi^+}=0$ and using (\ref{eq:ppr2:1}), we obtain:
$$
\tilde{\tilde\psi}  = \tilde\psi - \tilde f_2 
\frac{\omega_{\tilde\psi,\tilde f_2^+}}{\omega_{\tilde f_2,\tilde f_2^+}}= \psi - f_1 \frac{\omega_{\psi,f_1^+}}{\omega_{f_1,f_1^+}} - 
$$
$$
-\left(f_2 - f_1 \frac{\omega_{f_2,f_1^+}}{\omega_{f_1,f_1^+}}\right) 
\left(\frac{\omega_{\psi,f_2^+}\omega_{f_1,f_1^+}-\omega_{\psi,f_1^+}\omega_{f_1,f_2^+}}{\omega_{f_1,f_1^+}}\right)
\left(\frac{\omega_{f_2,f_2^+}\omega_{f_1,f_1^+}-\omega_{f_2,f_1^+}\omega_{f_1,f_2^+}}{\omega_{f_1,f_1^+}}\right)^{-1}=
$$
$$
= \psi  - f_1 \frac{1}{\omega_{f_1,f_1^+}\omega_{f_2,f_2^+} -\omega_{f_2,f_1^+}\omega_{f_1,f_2^+}} \times
$$
$$
\times\left(\omega_{\psi,f_1^+}\omega_{f_2,f_2^+} -
\frac{\omega_{\psi,f_1^+}\omega_{f_2,f_1^+}\omega_{f_1,f_2^+}}{\omega_{f_1,f_1^+}} 
-\omega_{\psi,f_2^+}\omega_{f_2,f_1^+}+
\frac{\omega_{f_2,f_1^+}\omega_{\psi,f_1^+}\omega_{f_1,f_2^+}}{\omega_{f_1,f_1^+}} \right)-
$$
$$
- f_2
\frac{1}{\omega_{f_1,f_1^+}\omega_{f_2,f_2^+} -\omega_{f_2,f_1^+}\omega_{f_1,f_2^+}} 
\left( \omega_{\psi,f_2^+}\omega_{f_1,f_1^+}-\omega_{\psi,f_1^+}\omega_{f_1,f_2^+} \right)=
$$
\begin{equation}
\label{eq:ppr2:2}
= \psi - \begin{bmatrix} f_1 & f_2  \end{bmatrix}
\frac{1}{\omega_{f_1,f_1^+}\omega_{f_2,f_2^+} -\omega_{f_2,f_1^+}\omega_{f_1,f_2^+}} 
\begin{bmatrix}
\omega_{\psi,f_1^+}\omega_{f_2,f_2^+} - \omega_{\psi,f_2^+}\omega_{f_2,f_1^+} \\
\omega_{\psi,f_2^+}\omega_{f_1,f_1^+} - \omega_{\psi,f_1^+}\omega_{f_1,f_2^+}
\end{bmatrix}.
\end{equation}
Taking into account that
$$
\Omega=\begin{bmatrix} 
\omega_{f_1,f_1^+} & \omega_{f_2,f_1^+}  \\
\omega_{f_1,f_2^+} & \omega_{f_2,f_2^+} 
\end{bmatrix}, \  
\Omega^{-1}=\frac{1}{\omega_{f_1,f_1^+}\omega_{f_2,f_2^+} -\omega_{f_2,f_1^+}\omega_{f_1,f_2^+}} 
\begin{bmatrix} 
\omega_{f_2,f_2^+} & -\omega_{f_2,f_1^+}  \\
-\omega_{f_1,f_2^+} & \omega_{f_1,f_1^+} 
\end{bmatrix},
$$
$$
\Omega^{-1}\begin{bmatrix}\omega_{\psi,f_1^+} \\ \omega_{\psi,f_2^+} \end{bmatrix} =
\frac{1}{\omega_{f_1,f_1^+}\omega_{f_2,f_2^+} -\omega_{f_2,f_1^+}\omega_{f_1,f_2^+}} 
\begin{bmatrix}
\omega_{\psi,f_1^+} \omega_{f_2,f_2^+}- \omega_{\psi,f_2^+} \omega_{f_2,f_1^+} \\
\omega_{\psi,f_2^+}\omega_{f_1,f_1^+}- \omega_{\psi,f_1^+}\omega_{f_1,f_2^+}
\end{bmatrix},
$$
formula (\ref{eq:ppr2:2}) can be rewritten as:
\begin{equation}
\label{eq:ppr2:3}
\tilde{\tilde\psi} = \psi  - \begin{bmatrix}f_1 & f_2 \end{bmatrix} 
\Omega^{-1}\begin{bmatrix}\omega_{\psi,f_1^+} \\ \omega_{\psi,f_2^+} \end{bmatrix}. 
\end{equation}
One can see that (\ref{eq:ppr2:3}) coincides with formula (\ref{eq:moutard5}) for $N=2$. 

The computations for $\tilde{\tilde\psi}^+$ are similar.

In additions, due to the definition of $\mathcal{M}$ according to (\ref{eq:pr2:moutard1})-(\ref{eq:comp1}),
using (\ref{eq:trans:pot}) for $\omega_{\tilde f_2,\tilde f_2^+}$  with $c_{\tilde f_2,\tilde f_2^+}=0$ 
and using (\ref{eq:ppr2:1}), we obtain:
$$
\tilde{\tilde u}= \tilde u + \frac{\tilde f_2 \overline{\tilde f_2^+}}{\omega_{\tilde f_2,\tilde f_2^+}} = 
u + \frac{f_1\overline{f_1^+}}{\omega_{f_1,f_1^+}} + 
\frac{\tilde f_2 \overline{\tilde f_2^+}}{\omega_{\tilde f_2,\tilde f_2^+}} =
$$
$$
=u + \frac{f_1\overline{f_1^+}}{\omega_{f_1,f_1^+}} + 
\left(f_2 -f_1 \frac{\omega_{f_2,f_1^+}}{\omega_{f_1,f_1^+}} \right) 
\left(\overline{f_2^+} -\overline{f_1^+} \frac{\omega_{f_1,f_2^+}}{\omega_{f_1,f_1^+}} \right) 
\frac{\omega_{f_1,f_1^+}}{\omega_{f_1,f_1^+}\omega_{f_2,f_2^+}-\omega_{f_2,f_1^+}\omega_{f_1,f_2^+}}= 
$$
$$
=u + \frac{ f_1\overline{f_1^+}\, \omega_{f_2,f_2^+} -f_2\overline{f_1^+}\, \omega_{f_1,f_2^+} -  
f_1\overline{f_2^+}\, \omega_{f_2,f_1^+} + f_2\overline{f_2^+}\, \omega_{f_1,f_1^+} }
{\omega_{f_1,f_1^+}\omega_{f_2,f_2^+}-\omega_{f_2,f_1^+}\omega_{f_1,f_2^+}}=
$$
\begin{equation}
\label{eq:ppr2:4}
=u+ \frac{1}{\omega_{f_1,f_1^+}\omega_{f_2,f_2^+}-\omega_{f_2,f_1^+}\omega_{f_1,f_2^+}}
\begin{bmatrix} f_1 & f_2 \end{bmatrix} \begin{bmatrix} 
\omega_{f_2,f_2^+} & -\omega_{f_2,f_1^+}  \\
-\omega_{f_1,f_2^+} & \omega_{f_1,f_1^+} 
\end{bmatrix}
\begin{bmatrix}\overline{f_1^+} \\\overline{f_2^+} \end{bmatrix}.
\end{equation}
One can see that (\ref{eq:ppr2:4}) coincides with formula (\ref{eq:moutard3}) for $N=2$.

This completes the proof of Theorem~\ref{th:1}.

\subsection{Proof of  Thereom~\ref{th:2}}
First, we check that $\omega_{\tilde f_2,\tilde\psi^+}$, $\omega_{\tilde\psi,\tilde f_2^+}$, $\omega_{\tilde f_2,\tilde f_2^+}$ defined in 
(\ref{eq:s3:moutard5}) are the  
potentials for the pairs  $\{{\hat{f}}$, ${\tilde\psi}^+\}$, 
$\{{\tilde\psi}$, ${\hat{f}}^+\}$, $\{{\hat{f}}$, ${\hat{f}}^+\}$:
$$
\partial_z \omega_{\tilde\psi,\tilde f_2^+}=-i\partial_z \frac{\omega_{\psi,f_1^+}}{\omega_{f_1,f_1^+}}=
-i\frac{\psi f_1^+}{\omega_{f_1,f_1^+}}+i\,\frac{\omega_{\psi,f_1^+}}{\omega^2_{f_1,f_1^+}}\,
f_1f_1^+=
$$
$$
=\left(\psi-f_1
\frac{\omega_{\psi,f_1^+}}{\omega_{f_1,f_1^+}}\right)\left( \frac{-if_1^+}{\omega_{f_1,f_1^+}}\right)=
\tilde\psi\, \hat{f}^+,
$$
$$
\partial_{\bar z} \omega_{\tilde\psi,\tilde f_2^+}=-i\partial_{\bar z} \frac{\omega_{\psi,f_1^+}}{\omega_{f_1,f_1^+}}=
i\frac{\overline{\psi}\overline{f_1^+}}{\omega_{f_1,f_1^+}}-i\,\frac{\omega_{\psi,f_1^+}}{\omega^2_{f_1,f_1^+}}\,
\overline{f_1}\overline{f_1^+}=
$$
\begin{equation}
\label{eq:ppr3:1}
=-\left(\overline{\psi}-\overline{f_1}
\frac{\omega_{\psi,f_1^+}}{\omega_{f_1,f_1^+}}\right)\left( \frac{-i\overline{f_1^+}}{\omega_{f_1,f_1^+}}\right)=
-\overline{\tilde\psi}\, \overline{\hat{f}^+};
\end{equation}
the calcualtions for $\omega_{\tilde f_2,\tilde\psi^+}$ are analogous to the calculations for 
 $\omega_{\tilde \psi,\tilde f_2^+}$ ;
$$ 
\partial_{z}\omega_{\tilde f_2,\tilde f_2^+}=\partial_{z}\frac{1}{\omega_{f_1,f_1^+}}= 
-\frac{\partial_{z}\omega_{f_1,f_1^+}}{\omega^2_{f_1,f_1^+}}=
-\frac{f_1f_1^+}{\omega^2_{f_1,f_1^+}}= 
\hat{f}\, \hat{f}^+,
$$
\begin{equation}
\label{eq:ppr3:2}
\partial_{\bar z}\omega_{\tilde f_2,\tilde f_2^+}=\partial_{\bar z}\frac{1}{\omega_{f_1,f_1^+}}= 
-\frac{\partial_{\bar z}\omega_{f_1,f_1^+}}{\omega^2_{f_1,f_1^+}}=
\frac{\overline{f_1}\,\overline{f_1^+}}{\omega^2_{f_1,f_1^+}}
=-\overline{\hat f}\,
\overline{\hat \beta^+}.
\end{equation}
Here, we used also that all potentials $\omega_{\psi,\psi^+}$ are pure imaginary.

Second, we calculate the transform $\mathcal{M}_2 \circ \mathcal{M}_1$:
$$
\tilde{\tilde\psi} = \tilde\psi -
\hat{f}\, \frac{\omega_{\tilde \psi,\tilde f_2^+}}{\omega_{\tilde f_2,\tilde f_2^+}}=\psi - f_1\,
\frac{\omega_{\psi,f_1^+}}{\omega_{f_1,f_1^+}}- \left(-i\frac{f_1}{\omega_{f_1,f_1^+}}\right) 
\frac{-i\,\frac{\omega_{\psi,f_1^+}}{\omega_{f_1,f_1^+}}}{\frac{1}{\omega_{f_1,f_1^+}}}=
$$
$$
= \psi  - f_1 \frac{\omega_{\psi,f_1^+}}{\omega_{f_1,f_1^+}}+ f_1 \frac{\omega_{\psi,f_1^+}}{\omega_{f_1,f_1^+}}= 
\psi;
$$
the calcualtions for $\tilde{\tilde\psi}^+$ are analogous to the calcualations for 
$\tilde{\tilde\psi}$;
$$
\tilde{\tilde u}= \tilde u + \frac{\hat{f} \overline{\hat{f}^+}}{\omega_{\tilde f_2,\tilde f_2^+}} = 
u + \frac{f_1\overline{f_1^+}}{\omega_{f_1,f_1^+}} + 
\frac{\hat{f} \overline{\hat{f}^+}}{\omega_{\tilde f_2,\tilde f_2^+}} =
u + \frac{f_1\overline{f_1^+}}{\omega_{f_1,f_1^+}} + 
\frac{-i\frac{f_1}{\omega_{f_1,f_1^+}} \overline{\left[-i\frac{f_1^+}{\omega_{f_1,f_1^+}}\right]} }
{\frac{1}{\omega_{f_1,f_1^+}}}=
$$
\begin{equation}
\label{eq:ppr3:3}
=u + \frac{f_1\overline{f_1^+}}{\omega_{f_1,f_1^+}} - 
\frac{f_1 \overline{f_1^+}}{\omega_{f_1,f_1^+}} = u.
\end{equation}

This completes the proof of Thereom~\ref{th:2}.

\section{Proof of meromorphic class conditions}
\label{sec:5}

We consider equation (\ref{eq:conf4}) and formulas (\ref{eq:sing2}), (\ref{eq:sing3}) for $u_*$, $\psi_*$,
where $n=1$, $r_{-1}=-1/2$, $n'=1$. In this case formulas  (\ref{eq:sing2}), (\ref{eq:sing3}) can be written as:
\begin{equation}
\label{eq:sec5:1}
u_*=e^{2i\phi(y)}\left(-\frac{1}{2x}+r_o(y)+r_1(y) x+ r_2(y) x^2+\ldots   \right),
\end{equation}
\begin{equation}
\label{eq:sec5:2}
\psi_*=e^{i\phi(y)}\left(\frac{\beta_{-1}(y)}{x}+\beta_0(y)+\beta_1(y) x+ \beta_2(y) x^2+\ldots   \right).
\end{equation}

We substitute (\ref{eq:sec5:1}), (\ref{eq:sec5:2}) into  (\ref{eq:conf4}), 
and we use that 
\begin{align}
\label{eq:sec5:3}
2\partial_{\bar\tau}\psi_*&=e^{i\phi(y)}\left(-\frac{\beta_{-1}(y)}{x^2}+
\frac{i \beta'_{-1}(y)-\phi'(y)\beta_{-1}(y) }{x}+ 
\vphantom{\sum\limits_{k=0}^{+\infty}} \right.\\ 
&+\left.\sum\limits_{k=0}^{+\infty}\left[ i \beta'_{k}(y) -\phi'(y)\beta_{k}(y) + 
(k+1)\beta_{k+1}(y) \right]x^k \right),\nonumber
\end{align}
\begin{align}
\label{eq:sec5:4}
2 u_*\psi_* &=e^{i\phi(y)}\left(-\frac{\overline{\beta_{-1}(y)}}{x^2}+
\frac{-\overline{\beta_{0}(y)} + 2r_0(y)\overline{\beta_{-1}(y)} }{x}+
\vphantom{\sum\limits_{k=0}^{+\infty}} \right.\\ 
&+\left.\sum\limits_{k=0}^{+\infty} \left[- \overline{\beta_{k+1}(y)} 
+ 2r_{k+1}(y)\overline{\beta_{-1}(y)} 
+2\sum\limits_{l=0}^{k}  r_l(y)\overline{\beta_{k-l}(y)} \right]x^k \right).
\nonumber
\end{align}
From this point and till the end of the proof $\, '$ denotes $\partial_y$.

Collecting the terms at $x^k$, $k=-2,-1,0,1,2,\ldots$, we obtain:
\begin{align}
\label{eq:sec5:5}
\overline{\beta_{-1}(y)}=&\beta_{-1}(y)& &\ \mbox{for} \  k=-2,\\
\label{eq:sec5:6}
\overline{\beta_{0}(y)}=&-i \beta'_{-1}(y) +\phi'(y)\beta_{-1}(y) + 
2r_0(y)\overline{\beta_{-1}(y)}& &\ \mbox{for} \ k=-1,\\
\label{eq:sec5:7}
\overline{\beta_{1}(y)}+\beta_{1}(y)=&-i \beta'_{0}(y) +\phi'(y)\beta_{0}(y) +
2r_1(y)\overline{\beta_{-1}(y)}+  & & \\
\nonumber  
& + 2r_0(y)\overline{\beta_{0}(y)} & & \ \mbox{for} \ k=0,\\
\label{eq:sec5:8}
\overline{\beta_{k+1}(y)}+(k&+1)\beta_{k+1}(y)=-i \beta'_{k}(y) +
\phi'(y)\beta_{k}(y) + & &\\
\nonumber
+& 2r_{k+1}(y)\overline{\beta_{-1}(y)} 
+2\sum\limits_{l=0}^{k}  r_l(y)\overline{\beta_{k-l}(y)} & & \ \mbox{for} \ k\ge1.
\end{align}

One can see that: relation (\ref{eq:sec5:5}) coincides with (\ref{eq:sing7});
relation (\ref{eq:sec5:6}) can be considered as a formula for finding $\beta_0$;
relations (\ref{eq:sec5:8}) can be considered as recursion relations for finding
$\beta_j$, $j\ge2$. In addition, the real part of (\ref{eq:sec5:7}) can be 
considered as a formula for finding $\Re\beta_{1}$, whereas the imaginary part
of (\ref{eq:sec5:7}) can be rewritten as 
\begin{equation}
\label{eq:sec5:9}
\Im\left[-i \beta'_{0}(y) +
\phi'(y)\beta_{0}(y)+ 2r_{1}(y)\overline{\beta_{-1}(y)} 
+2 r_{0}(y)\overline{\beta_{0}(y)}\right]=0.
\end{equation}
Actually, relations (\ref{eq:sec5:5}), (\ref{eq:sec5:6}),  (\ref{eq:sec5:9}) 
yield the solvability constraints on $\phi$, $r_0$, $r_1$. Substituting 
(\ref{eq:sec5:5}), (\ref{eq:sec5:6}) into (\ref{eq:sec5:9}) we obtain:
\begin{align}
\nonumber
&\Im\Big[i \phi'(y) \beta'_{-1}(y) + ((\phi'(y))^2 + 2 \phi'(y) \overline{r_0(y)})\beta_{-1}(y) 
+ 2r_1(y)\beta_{-1}(y)-\\ 
\label{eq:sec5:10}
& -2i r_0(y) \beta'_{-1}(y)
 +(2r_0(y)\phi'(y)+ 4(r_0(y))^2)\beta_{-1}(y) + \beta''_{-1}(y) -\\
\nonumber
&-i (\phi'(y) + 2\overline{r_0(y)})\beta'_{-1}(y) -
i (\phi''(y)+ 2\overline{r'_0(y)})\beta_{-1}(y)\Big]=0.
\end{align}
In turn, (\ref{eq:sec5:10}) can be rewritten as:
\begin{align}
\label{eq:sec5:11}
&\Im\Big[\beta''_{-1}(y) -  2i[ r_0(y) + \overline{r_0(y)}]  \beta'_{-1}(y)+
[(\phi'(y))^2+ \\
\nonumber
& + 2 \phi'(y) (\overline{r_0(y)} +r_0(y) )+ 
4(r_0(y))^2 + 2r_1(y)  -i \phi''(y) -2i \overline{r'_0(y)} ]\beta_{-1}(y) 
\Big]=0.
\end{align}
In addition, taking into account that $\phi$, $\beta_{-1}$ are real-valued, we
simplify (\ref{eq:sec5:11}) as follows:
\begin{equation}
\label{eq:sec5:12}
\Im\Big[-  2i[ r_0(y) + \overline{r_0(y)}]  \beta'_{-1}(y)+
[ 4(r_0(y))^2 + 2r_1(y)  -i \phi''(y) -2i \overline{r'_0(y)} ]
\beta_{-1}(y) \Big]=0.
\end{equation}
One can see that (\ref{eq:sec5:12}) is fulfilled for all sufficiently regular 
real-valued $\beta_{-1}$ if and only if (\ref{eq:sing5}), (\ref{eq:sing6}) are
fulfilled.

Finally, using (\ref{eq:sec5:5})-(\ref{eq:sec5:8}) under conditions 
(\ref{eq:sing5}), (\ref{eq:sing6}), one can see that all $\psi_*$ of 
(\ref{eq:sec5:2}) satisfying (\ref{eq:conf4}) with $u_*$ of \ref{eq:sec5:1})
are parametrised  by $\beta_{-1}$ and $\Im\beta_{1}$.
 
This completes the proof.

\section{Proof of Thereom~\ref{th:3}}
\label{sec:6}

Substituting (\ref{eq:sing7.5}), (\ref{eq:sing9}), (\ref{eq:sing10}) into (\ref{eq:conf4}), (\ref{eq:conf5}) 
we obtain:
\begin{equation}
\label{eq:sec6:1}
\beta_{0}(y)= i (\beta_{-1}(y))' + (\phi'(y) - 2 r_0(y) )\beta_{-1}(y), \ \ y\in[a,b],
\end{equation}
\begin{equation}
\label{eq:sec6:2}
\beta^+_{0}(y)= i (\beta^+_{-1}(y))' + (-\phi'(y) + 2 r_0(y))\beta^+_{-1}(y),\ \ y\in[a,b].
\end{equation}
As in Section~\ref{sec:5} we assume that $'$ denotes $\partial_y$.

Using  (\ref{eq:sing9}), (\ref{eq:sing10}) and (\ref{eq:sec6:1}), (\ref{eq:sec6:2}) we obtain:
\begin{align}
\nonumber
f_*f_*^+=&e^{-i\pi/2}\left(\frac{\beta_{-1}(y) \beta^+_{-1}(y) }{x^2}+
\frac{\beta_0(y) \beta^+_{-1}(y) + \beta_{-1}(y) \beta^+_{0}(y)  }{x} +O(1) \right)=\\
=& -i\frac{\beta_{-1}(y) \beta^+_{-1}(y)}{x^2}+\frac{1}{x}\Big( \beta_{-1}(y) \beta^+_{-1}(y) \Big)'+O(1)  
\ \ \mbox{in} \ \ \cT\cup\partial\cT.
\label{eq:sec6:3}
\end{align}

Next, equation (\ref{eq:kern1}) can be rewritten as:
\begin{equation}
\label{eq:sec6:4}
\partial_x\omega_{\psi,\psi^+}=2i\Im\Big( \psi\psi^+ \Big), \ \ 
\partial_y\omega_{\psi,\psi^+}=2i\Re\Big(\psi\psi^+ \Big). 
\end{equation}
As a corollary of (\ref{eq:sec6:3}),  (\ref{eq:sec6:4}) and the 
property that $\beta_1$,  $\beta_1^+$ are real-valued, we have:
\begin{equation}
\label{eq:sec6:5}
\partial_x\omega_{f_*,f_*^+}=\frac{-2i\beta_{-1}(y) \beta^+_{-1}(y)}{x^2}+ O(1), \ 
\partial_y\omega_{f_*,f_*^+}=\frac{2i\big(\beta_{-1}(y) \beta^+_{-1}(y)\big)'}{x}+ O(1);
\end{equation}
\begin{equation}
\label{eq:sec6:6}
\omega_{f_*,f_*^+}=\frac{2i}{x}\beta_{-1}(y) \beta^+_{-1}(y)+O(1) \ \ \mbox{in} \ \ \cT\cup\partial\cT.
\end{equation}
\begin{remark}
Note that
$$
\res\Big|_{x=0} \partial_x(\omega_{f_*,f_*^+})=0.
$$
for any fixed $y=y_0$.
\end{remark}

Using also the strict positivity of $\beta_{-1} \beta^+_{-1}$ we obtain: 
\begin{equation}
\label{eq:sec6:7}
\frac{1}{\omega_{f_*,f_*^+}}=-\frac{ix}{2\beta_{-1}(y) \beta^+_{-1}(y)}+O\left(x^2\right)  
\ \ \mbox{in} \ \ \cT_{a,b,\delta}\cup\partial\cT_{a,b,\delta},
\end{equation}
for some $\delta\in]0,\varepsilon[$.

Note also that
\begin{equation}
\label{eq:sec6:8}
f_*\overline{f_*^+}=i e^{2i\phi(y)} \frac{\beta_{-1}(y) \beta^+_{-1}(y)}{x^2}+ O\left(\frac{1}{x}\right) 
\ \ \mbox{in} \ \ \cT\cup\partial\cT.
\end{equation}
Finally, due to (\ref{eq:s2:moutard3}), (\ref{eq:sing7.5}), (\ref{eq:sec6:7}), (\ref{eq:sec6:8}):   
\begin{align}
\label{eq:sec6:9}
\nonumber
\tilde u_* &= u_* + \frac{f_*\overline{f_*^+}}{\omega_{f_*,f_*^+}}=e^{2i\phi(y)}\left(-\frac{1}{2x}+r_0(y)+ O(x) \right) + e^{2i\phi(y)}\frac{1}{2x} + O(1)= \\
&=O(1) \ \ \mbox{in} \ \ \cT_{a,b,\delta}\cup\partial\cT_{a,b,\delta}.
\end{align}
Theorem~\ref{th:3} is proved.


\begin{thebibliography}{ccc}

\bibitem {APP} V.A. Arkad'ev,  A.K. Pogrebkov,  M.K. Polivanov, Singular solutions of the
KdV equation and the inverse scattering method, \textit{Journal of Soviet Mathematics}, 
\textbf{31}(6) (1985), 3264-3279; doi:10.1007/BF02107228.

\bibitem{BC} R. Beals, R.R. Coifman, The spectral problem for the Davey-Stewartson and 
Ishimori hierarchies, In: ``Nonlinear evolution equations: integrability and spectral methods'', 
Proc. Workshop, Como/Italy 1988, Proc. Nonlinear Sci., (1990), pp. 15-23.
 
\bibitem{Bers} L. Bers, {\it Theory of pseudo-analytic functions},
Courant Institute of Mathematical Sciences, New York University,
Institute for Mathematics and Mechanics, 1953, 187 pages.

\bibitem{Crum} M.M. Crum, Associated Sturm-Liouville systems, \textit{Quart. J. Math. Oxford Ser. (2)} 
\textbf{6} (1955), 121-127; doi:10.1093/qmath/6.1.121.

\bibitem{DGNS} A. Doliwa, P. Grinevich, M. Nieszporski, P.M. Santini, Integrable lattices and their 
sublattices: From the discrete Moutard (discrete Cauchy-Riemann) 4-point equation to the self-adjoint 
5-point scheme, \textit{ J. Math. Phys.}, \textbf{48}(1) (2007), 013513; doi:10.1063/1.2406056.

\bibitem{DG}  J.J. Duistermaat,  F.A. Gr\"unbaum,  Differential equations in the spectral parameter,
\textit{Comm. Math. Phys.}, \textbf{103}(2) (1986),  177-240; doi:10.1007/BF01206937 .

\bibitem{F1} L.D. Faddeev, Inverse problem of quantum scattering theory. II,
\textit{Journal of Soviet Mathematics}, {\bf 5}(3) (1976), 334--396.

\bibitem{FS}  A.S. Fokas, L.-Y. Sung, On the solvability of the N-wave, Davey-Stewartson and 
Kadomtsev-Petviashvili equations, \textit{Inverse Problems}, \textbf{8}(5) (1992), 673-708;
doi:10.1088/0266-5611/8/5/002.

\bibitem{Gr} P.G. Grinevich, The scattering transform for the two-dimensional 
Schr\"odinger operator with a potential that decreases at infinity at 
fixed nonzero energy, \textit{Russian Math. Surveys}, \textbf{55}(6) (2000), 1015-1083; 
doi:10.1070/RM2000v055n06ABEH000333. 

\bibitem{GRN} P.G. Grinevich, R.G. Novikov, Transparent potentials at fixed energy in dimension 
two. Fixed-energy dispersion relations for the fast decaying potentials, 
\textit{Commun. Math. Phys.}, \textbf{174} (1995), 409-446; doi:10.1007/BF02099609.

\bibitem{GRN1}P.G. Grinevich, R.G. Novikov, Faddeev eigenfunctions for point potentials in 
two dimensions, \textit{Physics Letters A}, \textbf{376} (2012), 1102-1106; doi:10.1016/j.physleta.2012.02.025 .

\bibitem{GRN2}P.G. Grinevich, R.G. Novikov, Faddeev eigenfunctions for multipoint potentials, 
\textit{Eurasian Journal of Mathematical and Computer Applications}, \textbf{1}(2) (2013), 76-91.

\bibitem{GRN3} P.G. Grinevich, R.G. Novikov, Moutard transform for generalized analytic functions,
\textit{Journal of Geometric Analysis}, doi:10.1007/s12220-015-9657-8,  arXiv:1510.08764.

\bibitem{GRN4} P.G. Grinevich, R.G. Novikov, Generalized analytic functions, Moutard-type transforms and 
holomorphic maps, arXiv:1512.00343.

\bibitem{GN} P.G. Grinevich, S.P. Novikov, Two-dimensional ``inverse scattering problem'' for negative energies 
and generalized-analytic functions. 1. Energies below the ground state, 
\textit{ Functional Analysis and Its Applications}, \textbf{22}(1) (1988), 19-27; doi:10.1007/BF01077719. 

\bibitem{GN1} P.G. Grinevich, S.P. Novikov, Singular soliton operators and indefinite metrics,
\textit{Bulletin of the Brazilian Mathematical Society, New Series}, \textbf{44}(4)  (2013), 809-840; 
doi:10.1007/s00574-013-0035-5.

\bibitem{GN2} P.G. Grinevich, S.P. Novikov, Spectrally meromorphic operators and nonlinear systems,
\textit{Russian Mathematical Surveys}, \textbf{69}(5) (2014), 924-926; doi:10.1070/RM2014v069n05ABEH004922.

\bibitem{Kaz}  A.V. Kazeykina, A large-time asymptotics for the solution of the Cauchy problem 
for the Novikov-Veselov equation at negative energy with non-singular scattering data, 
\textit{Inverse Problems}, \textbf{28}(5) (2012), 055017, 21 pp; doi:10.1088/0266-5611/28/5/055017.

\bibitem{NKh} G.M. Henkin, R.G. Novikov,  The $\bar\partial$-equation in the 
multidimensional inverse scattering problem, \textit{Russian Math. Surveys} \textbf{42}(3)  (1987), 109-180;
doi:10.1070/RM1987v042n03ABEH001419.

\bibitem{Rnov2} E.L. Lakshtanov, R.G. Novikov, B.R. Vainberg, A global 
Riemann-Hilbert problem for two-dimensional inverse scattering at fixed energy,  arXiv:1509.06495.

\bibitem{MatvSal} V.B. Matveev, M.A. Salle, \textit{Darboux transformations and solitons}, 
Springer Series in Nonlinear Dynamics. Springer-Verlag, Berlin, 1991.

\bibitem{Mout} T.F. Moutard, Sur la construction des \'equations de la forme 
$\frac{1}{z}\frac{\partial^2 z}{\partial x\partial y}= \lambda (x, y)$ qui admettenent une int\'egrale 
g\'en\'erale explicite, \textit{J. \'Ecole Polytechnique}, \textbf{45} (1878), 1-11. 

\bibitem{NS}  J.J.C. Nimmo, W.K. Schief, Superposition principles associated with the Moutard 
transformation: an integrable discretization of a 2+1-dimensional sine-Gordon system, 
\textit{Proc. R. Soc. London A}, \textit{453} (1997), 255-279; doi:10.1098/rspa.1997.0015. 

\bibitem{RNov} R.G. Novikov, The inverse scattering problem on a fixed energy 
level for the two-dimensional Schr\"odinger operator, \textit{J. Funct. Anal.} 
\textbf{103}(2) (1992), 409-463; doi:10.1016/0022-1236(92)90127-5.

\bibitem{NT}  R.G. Novikov, I.A. Taimanov, The Moutard transformation and two-dimensional multipoint 
delta-type potentials, \textit{Russian Math. Surveys}, \textbf{68}(5) (2013), 957-959; 
doi:10.1070/RM2013v068n05ABEH004864

\bibitem{NTT} R.G. Novikov, I.A. Taimanov, S.P. Tsarev, 
Two-dimensional von Neumann-Wigner potentials with a multiple positive eigenvalue, 
\textit{Functional Analysis and Its Applications}, \textbf{48}(4) (2014), 295-297;
doi:10.1007/s10688-014-0073-9.

\bibitem{Taim1} I.A. Taimanov, Blowing up solutions of the modified Novikov-Veselov equation 
and minimal surfaces, \textit{Theoretical and Mathematical Physics}, \textbf{182}(2) (2015), 
173-181; doi:10.1007/s11232-015-0255-5.

\bibitem{Taim2}  I.A. Taimanov, The Moutard transformation of two-dimensional Dirac operators 
and M\"obius geometry, \textit{Mathematical Notes}, \textbf{97}(1) (2015), 124-135; 
doi:10.1134/S0001434615010149.

\bibitem{TT} I.A. Taimanov, S.P. Tsarev, On the Moutard transformation and its applications 
to spectral theory and Soliton equations, \textit{Journal of Mathematical Sciences}, 
\textbf{170}(3) (2010), 371-387; doi:10.1007/s10958-010-0092-x.

\bibitem{Vek} I.N. Vekua, {\it Generalized Analytic Functions}, Pergamon Press Ltd. 1962.

\end{thebibliography}
\end{document}